\documentclass{article}

\usepackage{amsfonts, amsmath, amssymb, amscd} \newfont{\eufm}{eufm10}
\let\Box\square

\begin{document}

\title{Algebraic entropy of elementary amenable groups.}  \author{D.V.  Osin
\thanks{The work has been supported by the RFFR grant 99-01-00894
and by the Swiss National Science Foundation.}} \maketitle
\begin{abstract}
We prove that any finitely generated elementary amenable group of
zero (algebraic) entropy contains a nilpotent subgroup of finite
index or, equivalently, any finitely generated elementary amenable
group of exponential growth is of uniformly exponential growth. We
also show that $0$ is an accumulation point of the set of
entropies of elementary amenable groups.
\end{abstract}


\section{Introduction.}


Let $G$ be a group generated by a finite set $X$. As usual, we
denote by $||g||_X$ the {\it word length} of an element $g\in G$
with respect to $X$, i.e., the length of a shortest word over the
alphabet $X\cup X^{-1}$ which represents $g$.

In this paper we study the {\it growth function} $\gamma_G^X :
\mathbb N \longrightarrow \mathbb N$ of $G$ which is defined by
$$\gamma _G^X(n)=card\; \{ g\in G\; :\; ||g||_X\le n\} .$$ Originally
growth considerations in group theory were introduced in 50-th by
Efremovich \cite{Ef}, \v Svarc \cite{Sh}, and F\o lner \cite{Fol},
and (independently) in 60-th by Milnor \cite{JM} with motivations
from differential geometry and theory of invariant means.

The {\it exponential growth rate} of $G$ with respect to $X$ is
the number $$\omega (G,X) = \lim_{n \to \infty} \sqrt[n]{\gamma
_G^X(n)} .$$ The above limit exists by submultiplicativity of
$\gamma_G^X $  \cite[Theorem 4.9]{Walte}. The quantity $$\omega
(G)=\inf\limits_X \omega (G,X) $$ is called a {\it minimal
exponential growth rate} of $G$ (the infimum is taken over all
finite generating sets of $G $). Finally, the {\it (algebraic)
entropy } of the group $G$ is defined by the formula
$$h(G)=\log\omega(G).$$

This notion of entropy comes from geometry and should not be
confused with the notion of entropy for a pair $(G, \mu )$, where
$\mu $ is a symmetric probability measure on a group $G$, as
defined in \cite{Ave}. In particular, if $G$ is a fundamental
group of a compact Riemannian manifold of unit diameter, then
$h(G)$ is a lower bound for the topological entropy of the
geodesic flow of the manifold \cite{Manni}. The exponential growth
rates appear also in the study of random walks on the Cayley
graphs of finitely generated groups.  Details and backgrounds can
be found in \cite{GrH}, \cite{H}.

The group $G$ is said to be of {\it exponential growth} if $\omega
(G,X)>1$ and of {\it subexponential growth} if $\omega (G,X)=0$.
If there exist constants $C,d>0$ such that $\gamma _G^X(n)\le
Cn^d$ for all $n\in \mathbb N$, then $G$ is said to be {\it of
polynomial growth}. These definitions depend on the group $G$
only, not on the choice of finite generating sets. We refer to
\cite{Chou}, \cite{Gri-84}, \cite{M}, \cite{T}, and \cite{W}, for
classical results concerning growth of various classes of groups.
Further, one says that $G$ has {\it uniformly exponential growth}
if $\omega (G)>1$ (or, equivalently, $h(G)>0$). The following
important problem goes back to the book \cite{GLP} and can be
found in \cite{Gri-91} as well as in \cite{GrH} and \cite{H}.

{\bf Question 1.1.} {\it Does there exist a finitely generated
group of non--uniform exponential growth, i.e., of exponential
growth and of zero entropy?}

On one hand, the affirmative answer was recently obtained by J.
Wilson \cite{Wil}. On the other hand, there are many examples of
classes of groups which are known to have uniformly exponential
growth. Let us mention some of them.

$\bullet $ Hyperbolic groups containing no cyclic subgroups of
finite index \cite{Koubi}.

$\bullet $ Free products with amalgamations $G\ast _{A=B}H$
satisfying the condition $(|G:A|-1)(|H:B|-1)\ge 2$ and
HNN--extensions $G\ast _A$ associated with a monomorphism $\phi
_1, \phi _2: A\to G$, where $|G:\phi _1(A)| + |G:\phi _2 (A)|\ge
3$ (see \cite{BH}).

$\bullet $ One--relator groups of exponential growth
\cite{GrH1-rel}.

$\bullet $ Solvable groups of exponential growth \cite{Osin} (the
particular case of polycyclic groups was considered independently
in \cite{Alp}).

$\bullet $ Linear groups of exponential growth \cite{Mozes}.

The main goal of the present paper is to investigate the case of
elementary amenable groups and discuss certain applications of the
obtained results.

{\bf Acknowledgments. } I am grateful to Pierre de la Harpe for
for useful comments and for his hospitality at University of
Geneva where the present article has been written.

\section{Main results}

In order to explain the Hausdorff--Banach--Tarski paradox, von
Neumann \cite{Neu} introduced the class of amenable groups. He
showed that all finite and abelian groups are amenable and the
class of amenable groups,  $AG$, is closed under four standard
operations of constructing new groups from given ones:

(S) Taking subgroups.

(Q) Taking quotient groups.

(E) Group extensions.

(U) Direct limits (that is, for a given set of groups $\{ G_\lambda
\} _{\lambda \in \Lambda } $ such that for any $\lambda , \mu \in
\Lambda $, there is $\nu \in \Lambda $ satisfying $G_\lambda \cup
G_\mu \subseteq G_\nu$, one takes $\bigcup\limits_{\lambda \in
\Lambda } G_\lambda $).

As in \cite{Day}, let  $EG$ be the class of {\it elementary
amenable groups} that is the smallest class  which contains all
abelian and finite groups, and closed under the operations
(S)--(U). In particular, $EG$ contains all solvable groups.
However, it is easy to construct a finitely generated group $G\in
EG$ that is not even virtually solvable. (Recall that a group $G$
is a virtually $\mathcal P$ group, where $P$ is a class of groups,
if there exists a subgroup $H$ of finite index in $G$ such that
$H\in \mathcal P$.)

The main result of this paper is the following.

{\bf Theorem 2.1.} {\it Let $G$ be a finitely generated elementary
amenable group of zero entropy. Then $G$ contains a nilpotent
subgroup of finite index. In particular, any elementary amenable
group of exponential growth is of uniformly exponential growth.}

This extends the result of Chou, saying that any elementary amenable
group of subexponential growth contains a nilpotent subgroup of
finite index, as well as the result of the author from \cite{Osin},
where the analog of Theorem 2.1 was proved in the case of solvable
groups.

In \cite{R}, Rosset proved that if $G$ is a group of
subexponential growth, $H$ is a normal subgroup of $G$, and $G/H$
is solvable, then $H$ is finitely generated. The methods developed
in the present paper allows to obtain the following more general
result on the structure of normal subgroups of groups with zero
entropy.

{\bf Theorem 2.2.} {\it Let $G$ be a finitely generated group of
zero entropy and let $H$ be a subgroup of $G$ such that the
quotient group $G/H$ is elementary amenable. Then $H$ is finitely
generated.}

In particular, Theorem 2.2 provides a natural approach to prove
that a group has uniform exponential growth.

Given a closed manifold $M$ endowed with a Riemannian metric $g$,
we denote by $h_{top}(M,g)$ the topological entropy of the
geodesic flow on $M$ (precise definitions and certain properties
can be found in \cite{Manni}, \cite{AKM}). The connection between
the topological entropy of geodesic flows and homotopic properties
of Riemannian manifold was observed by Dinamburg in \cite{D},
where it was proved that $h_{top}(M, g)>0$ whenever $\pi _1(M)$
has exponential growth. It is an interesting question to describe
the behaviour of $h_{top}(M,g)$ for a given $M$, when the metric
varies. One can show that $h_{top}(M,g)$ can be made arbitrary
large by local variations of the metric. Thus the question above
is the question of the existence of a non--trivial lower bound for
the quantity
\begin{equation}
h(M)=\inf\limits_{Vol(M,g)=1} h_{top}(M,g), \label{minent}
\end{equation}
where the infimum is taken over all Riemannian metrics on $M$
normalized by the condition $Vol(M,g)=1$.

It is known that $h(M)>0$ for any manifold admitting negative
sectional curvature \cite{Kat}, \cite{BCG}. The proof of this
result given in \cite{BCG} is of purely geometric nature. Here we
illustrate another (algebraic) approach to the estimation of
$h(M)$ from below. According to the main result of \cite{Manni},
we have $$h(M)\ge h(\pi _1(M))$$ for every closed manifold $M$.
Therefore, Theorem 2.1 yields the following generalization of the
Dinamburg theorem in the particular case of manifolds with
elementary amenable fundamental groups.

{\bf Corollary 2.1.} {\it Let $M$ be a closed Riemannian manifold
and $h(M)$ denote the topological entropy of $M$ defined by
(\ref{minent}). Then $h(M)>0$ whenever the fundamental group $\pi
_1 (M) $ is elementary amenable and has exponential growth.}

We note that Corollary 2.1 provides purely topological conditions for
the positivity of the entropy. However, the inequality $h(M)>0$ can
be true even for simply--connected manifolds (see \cite{Bab}).

In conclusion we discuss a question concerning possible values of the
quantity $\omega (G)$, where $G$ is a finitely generated group. For a
class $\mathcal U$ of finitely generated groups, we set $$\Omega
(\mathcal U) =\{ \omega (G)\; : \; G\in \mathcal U\}.$$ Clearly
$\Omega (\mathcal U)\subseteq [1, \infty)$. It is an important
problem to describe the precise structure of the set $\Omega
(\mathcal U)$ for various classes $\mathcal U$. The main question we
are interested here is whether $1$ is an accumulation point of
$\Omega (\mathcal U)$. Certain considerations of such a kind can be
found in \cite{GrH}, where Grigorchuk and de la Harpe showed that
there exist finitely generated groups of exponential growth with
growth rates arbitrary close to $1$. However, there are numerous
questions, which are still open. Probably the most interesting one is

{\bf Question 2.1.} {\it Let $\mathcal H$ denotes the set of all
hyperbolic groups which are non-elementary (i.e., not
cyclic--by--finite). Is $1$ an accumulation point of $\Omega
(\mathcal U)$? }

It was pointed out in \cite{Osin1}, that the positive answer would
imply the existence of a relatively simple construction of a
finitely generated group of non--uniform exponential growth. We
note that the groups constructed in \cite{GrH} are not hyperbolic,
but they are non--amenable and semihyperbolic in the sense of
\cite{Bri}.

We obtain the following analog of the Grigorchuk -- de la Harpe
result in the case of elementary amenable groups.

{\bf Theorem 2.3} {\it The number $1$ is an accumulation point of
the set $\Omega (EG)$.}

In contrast, we provide one result of the converse type.

{\bf Theorem 2.4.} {\it Let $\mathcal M$ denote the set of all
finitely generated metabelian non--polycyclic groups. Then there
is $\varepsilon >0$ such that $\omega (G)>1+\varepsilon $ for any
$G\in \mathcal M$. }

It is worth to mention that we do not know whether the analog of
Theorem 2.4 is true in the case of polycyclic groups.


\section{Outline of the proof of main theorems.}


Here we describe shortly the main idea of the proof of Theorems
2.1 and 2.2. Recall that a group $P$ is called {\it polycyclic} if
there is a finite subnormal series
\begin{equation}
1=P_k\vartriangleleft P_{k-1}\vartriangleleft \ldots
\vartriangleleft P_0=P,\label{pol}
\end{equation}
where $P_{i-1}/P_i$ is cyclic for every $i=1, \ldots , k$. We
begin by proving the following 'dichotomy'.

{\bf Proposition 3.1.} {\it Let $G\in EG$ be a finitely generated
group. Suppose, in addition, that $G$ is not virtually nilpotent.
Then there exists a normal subgroup $H$ of $G$ such that $G/H$ is
virtually polycyclic and at least one of the following conditions
holds.

1) $G/H$ has exponential growth.

2) $H$ is not finitely generated.}

Thus the proof of Theorem 2.1 is divided into two parts depending
on the existence of a finite generating set of the subgroup $H$
from Proposition 3.1. The first case is relatively simple, however
we need two auxiliary results to treat it.

{\bf Lemma 3.2.} {\it Let $G$ be a finitely generated group. Then
the following assertions are true.

1) Suppose $R$ is a normal subgroup of $G$; then $\omega (G/R)\le
\omega (G)$.

2) Suppose $R$ is a subgroup of finite index in $G$; then $\omega
(R)\le \omega (G)^{(2[G:R]-1)}$.}

The proof of claim {\it 1)} is straightforward and is left as an
exercise. Claim {\it 2)} follows, for example, from Proposition
3.3 of \cite{SW} (and can also be proved by the reader by using
straightforward arguments).

{\bf Theorem 3.1.} \cite[Theorem 1.1.]{Osin} {\it Let $G$ be a
finitely generated solvable group of zero entropy. Then $G$
contains a nilpotent subgroup of finite index.}

By combining the second assertion of Lemma 3.1 and Theorem 3.1, we
obtain that any virtually polycyclic group of exponential growth
is of uniform exponential growth, since every finitely generated
virtually nilpotent group has polynomial growth \cite{Bass}.
Further, if a group $G$ has a quotient group of uniform
exponential growth, then it is of uniform exponential growth
itself by the first assertion of Lemma 3.1. This shows that $G$
has uniform exponential growth in case the condition 1) of
Proposition 3.1 holds.

The second case is more complicated. Without loss of generality,
we can assume the quotient group $G/H$ to be not of exponential
growth, i.e., to be virtually nilpotent. Moreover, since the
property to be of uniform exponential growth is preserved under
the taking of subgroups of finite index, we can assume $G/H$ to be
nilpotent. In this settings the following proposition plays the
crucial role in our proof.

{\bf Proposition 3.2.} {\it Let G be a finitely generated group
such that there exists an exact sequence $$ 1\longrightarrow
K\longrightarrow G\longrightarrow N\longrightarrow 1, $$ where $N$
is nilpotent of degree $d$ and $K$ is not finitely generated. Then
we have }
\begin{equation}
\omega (G)\ge \sqrt[\alpha ]{2},\label{24}
\end{equation}
where $\alpha =3\cdot 4^{d+1} $.

Now let us turn to Theorem 2.2. The conditions of the theorem
together with Theorem 2.1 and the first assertion of Lemma 4.1
imply that $G/H$ is virtually nilpotent. Therefore, $G$ contains a
subgroup $G_0$ of finite index such that $H\lhd G_0$ and $G_0/H$
is nilpotent. By the second assertion of Lemma 3.1, we have
$\omega (G_0)=1$. Therefore, $H$ is finitely generated by
Proposition 3.2

Thus to complete the proofs of Theorems 2.1 and 2.2 it remains to
prove Propositions 3.1 and 3.2. The first proposition is obtained
in the next section. The proof of the second one involves methods
based on commutator calculus. We give this proof in Section 6
modulo some auxiliary results, which are obtained in Section 5.


\section{Description and some properties of elementary amenable groups.}


First we recall the description of elementary classes of groups
given in \cite{Osin1}. In case of elementary amenable groups this
description is slightly stronger than the Chou one \cite{Chou}.

{\bf Definition 4.1.} Let $B$ be a class of groups. The {\it
elementary class of groups with the base} $B$ is the smallest
class of groups which contains $B$ and is closed under operations
(S)--(U).

Now we fix $B$. Let ${\mathcal E}_0(B)$ consist of the trivial
group only. Assume that $\alpha >0$ is an ordinal and that we have
defined ${\mathcal E}_\beta (B)$ for each ordinal $\beta < \alpha
$. If $\alpha $ is a limit ordinal, set $${\mathcal E}_\alpha
(B)=\bigcup\limits_{\beta <\alpha } {\mathcal E}_\beta (B),$$ and
if $\alpha  $ is successor, let ${\mathcal E}_\alpha (B)$ be the
class of groups that can be obtained from groups of ${\mathcal
E}_{\alpha -1}(B)$ by applying operation (U) or the following
operation once.

(E$_0$) Given a group, take its extension by a group from $B$.

{\bf Lemma 4.1.} \cite{Osin1} {\it The class $\mathcal E_\alpha $
is closed under operation (S) and (Q) (SQ--closed for brevity) for
each ordinal $\alpha $.}

{\bf Theorem 4.1.} \cite{Osin1} {\it Let $B$ be a class of groups.
Assume that $B$ is closed under the operations (S) and (Q). Then
we have $$ \mathcal E(B)=\bigcup\limits_{\alpha }\mathcal E_\alpha
(B), $$ where the union is taken over all ordinal numbers.}

{\bf Example 4.1.} Let us take $B=A \cup F$, where $A$ and $F$ are
the classes of all abelian and finite groups respectively. Then
the corresponding elementary class is precisely $EG$, as follows
from Theorem 4.1.

The next lemma is well--known (see for example \cite{Baer} or
\cite{Sch}).

{\bf Lemma 4.2.} {\it Any extension of a virtually polycyclic
group by a virtually polycyclic group is a virtually polycyclic
group.}

Now we are ready to formulate the main result of this section.
Instead of Proposition 3.1, we will prove a stronger result by
transfinite induction on $\alpha $. Recall that a subgroup $H$ of
a group $G$ is called characteristic if for any automorphism $\phi
$ of $G$, one has $\phi (H)\le H$. Evidently if $G$ is a normal
subgroup of a group $F$ and $H$ is a characteristic subgroup of
$G$, then $H$ is normal in $F$.

{\bf Lemma 4.3.} {\it Let $G\in EG_\alpha $ be a finitely
generated group. Suppose, in addition, that $G$ is not virtually
nilpotent. Then there exists a characteristic subgroup $H$ of $G$
such that $G/H$ is virtually polycyclic and at least one of the
following conditions holds.

1) $G/H$ has exponential growth.

2) $H$ is not finitely generated.}

{\it Proof.} The case $\alpha =0$ is trivial. Suppose that $\alpha
>0$. First assume that $\alpha $ is a limit ordinal. Then $G\in
EG_\beta $ for some $\beta <\alpha $ and thus the assertion of the
lemma holds by the inductive assumption.

Now let $\alpha $ be a non--limit ordinal. Assume that
$G=\bigcup\limits_{\lambda \in \Lambda } G_\lambda $, where
$G_\lambda \in EG_{\alpha -1}$. Since $G$ is finitely generated,
we have $G=G_{\lambda _0}$ for some $\lambda _0$. Hence the
assertion of the lemma is true by the inductive hypothesis again.
Further, suppose $G$ is an extension of the form $$
1\longrightarrow M\longrightarrow G\longrightarrow
L\longrightarrow 1, \label{ext} $$ where $M\in EG_{\alpha -1}$ and
$L$ is abelian or finite. First we consider the case of abelian
$L$. Take $G'=[G,G]$ and observe that $G'\le M$ and thus $G' \in
EG_{\alpha -1}$ by Lemma 4.1. If $G'$ is not finitely generated,
we can take it for $H$. Otherwise, there are two possibilities.
The first one is the case of virtually nilpotent $G'$. Clearly,
then $G$ is virtually polycyclic by Lemma 4.2. Taking into account
that $G$ is not virtually nilpotent, we conclude that $G$ is of
exponential growth.

Now suppose that $G'$ is not virtually nilpotent. Then, by the
inductive hypothesis, there exists a characteristic subgroup $H\le
G'$ satisfying the requirements of the proposition. Clearly, $H$
is a characteristic subgroup of $G$. It remains to notice that
since $G'/H$ is virtually polycyclic, then $G/H$  is virtually
polycyclic. Moreover, since $G'/H$ is of exponential growth, so is
$G/H$.

Similarly, if $L$ is finite, say $|L|=m$, then we take the
subgroup $G^m=\langle g^m\; : \; g\in G\rangle $. Evidently
$G^m\le M$ and the rest of the proof is essentially the same as in
the previous case. Additionally we only need the fact that any
finitely generated periodic group from $EG$ is finite (see
\cite{Chou}). $\Box $


\section{Technical lemmas.}


Throughout this section we fix a group $G$ generated by a finite
set $X$  and denote by ${\mathcal L}(X\cup X^{-1})$ the set of all
words over $X\cup X^{-1}$. For two words $u,v\in \mathcal L(X\cup
X^{-1})$, we write $u\equiv v$ to express the letter--for--letter
equality, and $u=v$ if $u$ and $v$ represent the same element of
$G$. By $\| w\| $ we denote the length of a word $w\in \mathcal
L(X\cup X^{-1})$. We write also $u^v$ instead of $v^{-1}uv$ and
$[u,v]$ instead of $u^{-1}v^{-1}uv$.

We fix arbitrary finite subsets $V, W\in {\mathcal L}(X\cup
X^{-1})$ and assume in addition that $V$ and $W$ satisfy the
following conditions.

(I) $X^{\pm 1}\in V$.

(II) The set $V$ is ordered, i.e., $V=\{ v_1, v_2, \ldots , v_p\}
$. Moreover, set  $V_i=\{ v_1, v_2, \ldots , v_i\} $ for each
$i=1,\ldots p $; then either $$[V_i, V_j]\subseteq V_{\min \{
i,j\} -1}$$ or $$[V_i, V_j]\subseteq W$$ for any $i,j=1,\ldots , p
$.

(III) For any $i=1,\ldots , p  $, and any $w\in W$, the normal
closure $$\langle w\rangle ^{\langle v_i\rangle }= \langle
v_i^{-l}wv_i^l\; : \; l\in \mathbb Z\rangle $$ is finitely
generated, i.e., there exists $L_i\in \mathbb N$ such that
$$\langle w\rangle ^{\langle v_i\rangle }= \langle
v_i^{-l}wv_i^l\; : \; |l|\le L_i \rangle .$$

We are going to show that $\langle W\rangle ^G$ is finitely
generated as a subgroup. To reach this this goal we need one more
auxiliary definition. Let $v$ be a word in the alphabet $V^{\pm
1}\cup W^{\pm 1}$. Denote by $\lambda _i(v)$ the number of
appearances of the letters $v_i^{\pm 1}$ in $v$. For instance, if
$v\equiv v_1v_2v_1^{-1}$, then $\lambda _1(v)=2,$ $ \lambda
_2(v)=1$. We note also that $\lambda _i $ is defined just for
words over $V^{\pm 1}\cup W^{\pm 1}$, not for elements of $G$, as
different words can represent the same element.

Given $X$, $G$, $V$, and $W$ as described above, we set $$
L=\max\limits_{i=1, \ldots , p} L_i $$ and
\begin{equation}
Z=\{ v^{-1}wv\; : \; w\in W,\; v\in\mathcal L(V^{\pm 1}),\;
\lambda _i (v)\le L\;  \forall i=1, \ldots , p\} . \label{Z}
\end{equation}
We will say that the above decomposition $z\equiv v^{-1}wv$ of a
word $z\in Z$, where $w\in W,\; v\in\mathcal L(V^{\pm 1})$, is a
{\it canonical form} of an element $z\in Z$; clearly $\lambda
_i(z)=2\lambda _i(v)$.

The main result of this section is the following.

{\bf Lemma 5.1.} {\it In the above notation, we have $\langle
W\rangle ^G=\langle Z\rangle $, i.e., the normal closure of $W$ in
the group $G$ is generated by the set $Z$ as a subgroup.}

The proof consists of four lemmas. Denote by $(x,y^{\pm 1})_1$ the
set $\{ [x,y], [x,y^{-1}]\} $. Further, we define $$ (x, y^{\pm
1})_{i+1}= \{ [c,y], [c, y^{-1}]\; :\; c\in (x, y^{\pm 1})_{i} \}.
$$

{\bf Lemma 5.2.} {\it Let $H$ be a group, $a,b\in H$. Then
$$(a,b^{\pm 1})_n\subseteq \langle a^{b^l}\;: \; l= -n, \ldots ,
n\} $$ for any $n\in \mathbb N$.}

{\it Proof.} For $n=1$, we have $$ (a, b^{\pm 1})_1= \{ a^{-1}a^b,
a^{-1}a^{b^{-1}}\} \subseteq \langle a, a^{b^{\pm 1}}\rangle . $$
Now suppose $n>1$. By induction, we can assume that the assertion
of the lemma is true for $(n-1)$, i.e.,
\begin{equation}
(a, b^{\pm 1})_{n-1}\subseteq \langle  a^{b^l}\; :\; l=-n+1,
\ldots , n-1 \rangle .\label{abn}
\end{equation}
Denote $a_l=a^{b^l}$ for brevity and consider an element
$$c=a_{l_1}^{\alpha _1}\ldots a_{l_m}^{\alpha _m}\in (a, b^{\pm
1})_{n-1}.$$ By (\ref{abn}), we can assume that $|l_j|\le n-1$ for
any $j=1, \ldots , m$. We obtain $$ [c, b]\equiv  c^{-1}c^b=
(a_{l_1}^{\alpha _1}\ldots a_{l_m}^{\alpha _m})^{-1}
(a_{l_1}^{\alpha _1}\ldots a_{l_m}^{\alpha _m})^b=(a_{l_1}^{\alpha
_1}\ldots a_{l_m}^{\alpha _m})^{-1} (a_{l_1+1}^{\alpha _1}\ldots
a_{l_m+1}^{\alpha _m}). $$ Therefore, $[c, b]\in \langle a_{-n+1},
\ldots , a_n\rangle $. Similarly we obtain $[c, b^{-1}]\in \langle
a_{-n}, \ldots , a_{n-1}\rangle $. The lemma is proved. $\Box $

The following three lemmas will be proved by common induction on a
parameter $r$.

{\bf Lemma 5.3.} {\it Let $0\le n_1< n_2<\ldots < n_m$ be a
sequence of integers. Consider a word
\begin{equation}
\bar v\equiv (a_1\ldots a_{n_1}v_r^{\epsilon _1})\cdot (
a_{n_1+1}\ldots a_{n_2}v_r^{\epsilon _2})\cdot \ldots \cdot
(a_{n_{m-1}+1}\ldots a_{n_m}v_r^{\epsilon _m}), \label{vvv}
\end{equation}
where $a_i\in (V\setminus \{ v_r\} )^{\pm 1}$, $\epsilon _i\in
\mathbb Z$ for each $i=1,\ldots , m$, and
\begin{equation}
\sum\limits_{i=1}^m|\epsilon _i|\le L+1. \label{eps}
\end{equation}
Then we have
\begin{equation}
\bar v=v_r^\sigma \cdot a_1b_1 \cdot a_2b_2\cdot \ldots \cdot
a_{n_m}b_{n_m}, \label{v}
\end{equation}
where $\sigma ={\sum\limits_{i=1}^m\epsilon _i}$ and
\begin{equation}
b_i\in \left\langle V_{r-1}\cup \left( \bigcup\limits_{j=-L}^L
W^{v_r^j}\right) \right\rangle  \label{bi}
\end{equation}
for all $i$. In particular, $b_i\in \langle V_{r-1}\cup Z\rangle
$.}

{\bf Lemma 5.4.} {\it Suppose that $z\equiv v^{-1}wv$ is a
canonical form of an element $z\in Z$ and $b$ is a word over
$(V\setminus \{ v_r\})^{\pm 1}\cup \left( \bigcup\limits_{j=-L}^L
W^{v_r^j}\right)^{\pm 1} $. Then we have
\begin{equation}
z^b \in \left\langle \begin{array}{lll} & & y_0\in \mathcal
L(V^{\pm 1}),\; w_0\in W,\\ y_0^{-1}w_0y_0 & : & \lambda _i(y_0)
\le L \; \forall \; i=1, \ldots , r,\\ & & \lambda _i(y_0)\le
\lambda _i(v)+\lambda _i(b)\; \forall \; i=r+1, \ldots , p
\end{array} \right\rangle . \label{za}
\end{equation}
In particular, if
\begin{equation}
\lambda _i(v)+\lambda _i(b)\le L \label{lambda}
\end{equation}
for each $i=r+1, \ldots, p$, then $z^b\in \langle Z\rangle .$}

{\bf Lemma 5.5.} {\it For any $z=v^{-1}wv\in Z$ and any $t\in
(V_r)^{\pm 1}$, we have
\begin{equation}
z^t \in \left\langle \begin{array}{lll}
& & y_0\in \mathcal L(V^{\pm 1}),\; w_0\in W,\\
y_0^{-1}w_0y_0 & : & \lambda _i(y_0) \le L \; \forall \; i=1, \ldots , r,\\
& & \lambda _i(y_0)\le \lambda _i(v)+\lambda _i(t)\; \forall \;
i=r+1, \ldots , p
\end{array} \right\rangle . \label{zt}
\end{equation}
In particular, $z^t\in \langle Z\rangle .$}

{\it Proof.} For all lemmas the case $r=1$ is essentially the same
as the inductive step. Thus we assume Lemmas 5.3 -- 5.5 to be true
for all positive integers $q<r$ whenever $r>1$ (and assume nothing
if $r=1$).

{\it Proof of Lemma 5.3.} Using the formula $$xy=yx[x,y]$$ we can
collect all appearances of the letter $v_r^{\pm 1}$ in the word
$\bar v$ from right to left in order to obtain a word of the form
(\ref{v}). It is easy to check that each $b_i$ will be a product
of elements of the sets $(a_i, v_r^{\pm 1})_n$, where
\begin{equation}
n\le \sum\limits_{i=1}^m |\epsilon _i| \label{nle}
\end{equation}
For an element $u\in (a_i, v_r^{\pm 1})_n$, there are two
possibilities.

(a) First assume that $u\in V$. Then $u\in V_{r-1}$ by condition
(II) (see the beginning of the section).

(b) Suppose $u\notin V$. Consider the minimal $n_0$ such that
$(a_i,v_r^{\pm 1})_{n_0}\not\subseteq V$. Clearly, $a_i\in V^{\pm
1}$ implies that $n_0\ge 1$. By Lemma A.2, $$ u\in \left( (a_i,
v_r^{\pm 1})_{n_0}, v_r^{\pm 1}\right)_{n-n_0}\subseteq
\left\langle \left( (a_i, v_r^{\pm 1})_{n_0}\right)^{v_r^l}\; : \;
|l|\le n-n_0\right\rangle . $$ Using (\ref{eps}) and (\ref{nle}),
we note that $$n-n_0\le \sum\limits_{i=1}^m |\epsilon _i| - n_0\le
L$$ and hence $$ u\in \left\langle \left( (a_i, v_r^{\pm 1})_{n_0}
\right)^{v_r^l}\; : \; |l|\le L\right\rangle \le \left\langle
\bigcup\limits_{j=-L}^LW^{v_r^j} \right\rangle . $$ Indeed, by
minimality of $n_0$, we have $(a_i, v_r^{\pm 1})_{n_0-1}\in V$.
Now using condition (II), we obtain $(a_i, v_r^{\pm 1})_{n_0}=
\left((a_i, v_r^{\pm 1})_{n_0-1}, v_r^{\pm 1}\right) \subseteq W$

Thus in both cases
$$
u\in \left\langle \bigcup\limits_{j=-L}^LW^{v_r^j} \cup
V_{r-1}\right\rangle ,
$$
and, therefore, the same is true for each $b_i$.  The lemma is
proved.

{\it Proof of Lemma 5.4.} Denote by $Y$ the group situating at the
right side of (\ref{za}). The proof will be by induction on the
length of the word $b$. The case $| b| =0$ is trivial. Now suppose
$| b| =n+1\ge 1$. Then $b=a_0a_1$, where $a_1\in (V\setminus \{
v_r\} )^{\pm 1}\cup \left( \bigcup\limits_{j=-L}^LW^{v_r^j}
\right) ^{\pm 1}$ and $a_0$ has length $n$. By the inductive
assumption, we have
\begin{equation}
z^b=z^{a_0a_1}=(z_1\ldots z_q)^{a_1}=z_1^{a_1}\ldots
z_q^{a_1},\label{zja1}
\end{equation}
where $z_j=y_j^{-1}w_jy_j$ are some elements such that $w_j\in
W^{\pm 1}$, $y_j\in \mathcal L(V^{\pm 1})$, $\lambda _i(y_j)\le
\lambda _i(a_0)+\lambda _i(v)$ for all $i=r+1, \ldots , p$, and
$\lambda _i(y_j)\le L$ for all $i=1, \ldots , r$.

Now let us consider $z_j^{a_1}$ for some $j$ and prove that
$z_j^{a_1}\in Y$. There are three possibilities.

(a) $a_1\in \left(\bigcup\limits_{j=-L}^L W^{v_r^j}\right) ^{\pm
1} $. Evidently $a_1\in Z$ in this case. Moreover, $a_1\in Y$ and
hence $z_j^{a_1}\in Y$.

(b) $a_1\in \left(V_{r-1}\right) ^{\pm 1}$. We note that this case
is impossible if $r=1$. If $r>1$, we assume that Lemma 5.5 has
already been proved for all smaller volumes of the parameter. Thus
we obtain $z_j^{a_1}\in Y$ applying Lemma 5.5 for $t\equiv a_1$,
$z\equiv z_j$.

(c) $a_1\in \left( V\setminus V_r\right) ^{\pm 1}$. Suppose
$a_1\equiv v_k$ for some $k\in \{r+1,\ldots , p\} $. For the
element $z_j^{a_1}$ consider its canonical form, the word
$v_k^{-1}y_j^{-1}w_jy_jv_k$, obtained from the canonical form of
the element $z_j$. We have $$ \lambda _k(v_ky_j)=\lambda _k(y_j)
+1\le \lambda _k(a_0)+\lambda _k(v)+1= \lambda _k(b)+\lambda
_k(v). $$ Clearly, if $i\ne k$, then $\lambda _i(v_ky_j)=\lambda
_i(y_j).$ This shows that $z_j^{a_1}$ lies in $Y$ again.

Since $z_j^{a_1}\in Y$ is true for each factor of type $z_j^{a_1}$
in (\ref{zja1}), we obtain $z^b\in Y$ and the proof of the lemma
is completed.

{\it  Proof of Lemma 5.5.} Denote by $F$ the group at the right
side of (\ref{zt}). In view of inductive arguments, it is
sufficient to consider the case $t\equiv v_r^{\pm 1}$. Assume that
$t\equiv v_r$ for convenience (the case $t\equiv v_r^{-1}$ is
analogous). First suppose that $\lambda _r(v)\le L-1,$ i.e.,
$\lambda _r(z)\le 2L-2.$ Note that $$ \lambda _i(z^{v_r})=\left\{
\begin{array}{l}
\lambda _i(z), \; {\rm if }\; i\ne r,\\
\lambda _i(z)+2, \; {\rm if }\; i=r.
\end{array}\right.
$$ Thus $\lambda _i (z^{v_r})\le 2L$ for $i=1, \ldots , r$, and
$\lambda _i(z^{v_r})= \lambda _i(z)$ for $i=r+1, \ldots ,p$. This
means that $z^{v_r}\in F$.

Now let $\lambda _r(v)=L.$ Then $\lambda _r (vv_r)=L+1$ and the
word $\bar v\equiv vv_r$ has the form (\ref{vvv}). Applying Lemma
5.3, we obtain $$\bar v=v_r^\sigma b,$$ where $$|\sigma
|=\left|\sum\limits_{i=1}^m \epsilon _i\right|\le \lambda _r(\bar
v)=L+1$$ and $b$ satisfies the condition $\lambda _i(b)=\lambda
_i(v)\le L$ for all $i=r+1, \ldots , p$ (obviously this condition
follows from (\ref{bi})).  In case $|\sigma |\le L$ we do nothing.
If $|\sigma |=L+1$, we apply condition (III) and obtain
\begin{equation}
z^{v_r}=w^{v_r^\sigma b}=\left( \prod\limits_{i=-L}^L \left(
w^{v_r^i}\right) ^{\xi _i}\right) ^b =\prod\limits_{i=-L}^L
\left(\left( w^{v_r^i}\right) ^b \right) ^{\xi _i}.
\end{equation}
Finally, we consider the elements $(w^{v_r^j})^b$, where $|j|\le
L$. In order to finish the proof of Lemma 5.5 it remains to show
that these elements belong to $F$. The element $b$ satisfies the
conditions of Lemma 5.4, as it contains no appearances of the
letters $v_r^{\pm 1}$. Thus $\left( w^{v_r^i}\right) ^b\in F$ by
Lemma 5.4. It follows that $z^{v_r}\in F$. The same arguments show
that $z^{v_r^{-1}}\in F$. The lemma is proved and the inductive
step is completed. $\Box $

{\it Proof of Lemma 5.1.} Lemma 5.5 implies that $z^t\in \langle
Z\rangle $ for any $z\in Z,$ $t\in V$. Since $X^{\pm 1}\subseteq
V$, we have $z^g\in \langle Z\rangle $ for any $g\in G$. This
means that $\langle Z\rangle ^G=\langle Z\rangle $ and the lemma
is proved. $\Box $


\section{Estimating exponential growth rates from below}


For any group $H$, let $\gamma _iH$ be the $i$-th term of the
lower central series $$H=\gamma _1H\vartriangleright \gamma
_2H\vartriangleright\ldots ,$$ where $\gamma _{i+1}H=[\gamma
_iH,H]$. Recall that a group $N$ is called {\it nilpotent of
degree $t$} if $\gamma _{t+1}N=1$. Finally, given subsets
$Y,Z\subseteq G$, let $\langle Y\rangle$ denote the subgroup
generated by $Y$, and $\langle Y \rangle ^Z$ the subgroup
generated by all elements of type $z^{-1}yz$, where $y\in Y, z\in
Z$. Thus $\langle Y\rangle ^G$ is the normal closure of $Y$ in
$G$.

{\bf Definition 6.1.} Let $G$ be a group with a given finite
generating set $X$. For any finite subset $Y=\{ y_1, \ldots ,
y_m\} \subseteq G$, we define its {\it depth} with respect to $X$
as follows $$ depth_X(Y)= \max\limits_{i=1,\ldots , m}||y_i||_X.
$$ If $H$ is a finitely generated subgroup of $G$, then we define
its depth with respect to $X$ by setting $$
depth_X(H)=\min\limits_{H=\langle Y\rangle } depth_X(Y), $$ where
the minimum is taken over all finite generating sets of $H$.

{\bf Lemma 6.1.} {\it Suppose that $G$ is a group with a given
finite generating set $X$ and $R$ is a finitely generated subgroup
of $G$; then we have $$\omega (G, X)\ge \left( \omega (R)\right)
^{\frac{1}{depth_X(R)} }.$$}

The proof is straightforward and is left as an exercise to the
reader. $\Box $

Let us introduce certain auxiliary notation. As above, suppose $G$ is
a group generated by a finite set $X$. Then we set $W_1(X)=X\cup
X^{-1}$ and, by induction,  $$W_i(X)=\{ [u^{\pm 1},v^{\pm 1}]\; :\;
u\in W(i_1), v\in W(i_2), i_1,i_2\in {\mathbb N}, i_1+i_2=i\} $$ for
any $i>1$. We write $weight\; (v)=i$ for a word $v\in \mathcal
L(X\cup X^{-1})$, if $v\in W_i$. Also, consider the function
$f:{\mathbb N}\to {\mathbb N}$ such that
\begin{equation}
f(1)=1 \;\;  {\rm and}\;\; f(n+1)=2f(n)+2  \label{fn}
\end{equation}
for any $n\in \mathbb N$. It can easily be checked that
$f(n)=3\cdot 2^{n-1}-2$. The following lemma is quite trivial.

{\bf Lemma 6.2.} {\it Let $f$ be the function given by (\ref{fn}).
Then for any $i,j\in \mathbb N$, one has $$2(f(i)+f(j))\le
f(i+j).$$}

{\bf Lemma 6.3.} {\it For any group $G$ with a given finite
generating set $X$, one has
\begin{equation}
depth _X(W_n(X))\le f(n).\label{depth}
\end{equation}}

{\it Proof.} We proceed by induction on $n$. The case $n=1$ is
trivial. Further, for $n>1$, we observe that if $u\in W_{i_1}(X),
v\in W_{i_2}(X)$ and $i_1+i_2=n$,  then $$\begin{array}{ll}
||[u^{\pm 1},v^{\pm 1}]||_X\le & 2(||u||_X + ||v||_X)\le 2(depth
_X(W_{i_1}(X))+depth _X(W_{i_2}(X))) \\ & \le 2(f(i_1)+f(i_2))\le
f(n)\end{array} $$ by the inductive hypothesis and Lemma 6.2.
$\Box
$

As an exercise, one can show that if $G$ is a non abelian free
group and $X$ is a basis in $G$, then $depth _X(W_n(X))=f(n)$.

The next lemma follows from the results of the previous section.

{\bf Lemma  6.4.} {\it Suppose that $G$ is a finitely generated
group and, for some $s\in\mathbb N$, $s\ge 2$, all subgroups of
type
\begin{equation}
H_{v,w}=\langle v^{-l}wv^l\; : \; l\in\mathbb Z\rangle
\label{Hv}
\end{equation}
are finitely generated for any $v\in \bigcup\limits_{j=1}^{s-1} W_j$,
$w\in\bigcup\limits_{j=s}^{2s} W_j$.
Then $\gamma _s(G) $ is finitely generated .}

{\it Proof.} The reader can easily check that the sets
$$W=\bigcup\limits_{j=s}^{2s} W_j$$ and
$$V=\bigcup\limits_{j=1}^{s-1} W_j$$ satisfy hypotheses (I) --
(III) listed at the beginning of Section 5. Indeed, (I) is
obvious. To satisfy (II), we just need to order commutators in $V$
in such a way that $weight (v_i)\le weight(v_j)$ whenever $i\ge
j$. Finally, (III) follows from the conditions of Lemma 6.4. It
remains to note that $\gamma _s(G)=\langle W\rangle ^G$, as
$\gamma _s(G)$ is generated by $\bigcup\limits_{t\ge s} W_t $ (see
\cite[Ch.5]{KMS}) and $W_t\in \langle W\rangle ^G$ for every $t>
2s$ by the definition of $W_t$. $\Box $

{\it Proof of Proposition 3.2.} Let $X$ be some finite generating
set of $G$. Let us put $s=d+1$. We would like to show that there
is a subgroup $H_{v,w}\le G$ of type (\ref{Hv}) having no finite
set of generators.  Indeed, suppose that all $H_{v,w}$ are
finitely generated. Then $\gamma _sG$  is finitely generated by
Lemma 6.4. Clearly $\gamma _sG \triangleleft K$. Therefore,
$K/\gamma _s G$ is a subgroup of a finitely generated nilpotent
group $G/\gamma _s G$ and thus is finitely generated. It follows
that $K$ is finitely generated and we arrive at contradiction.

Thus there exists $H_{v,w}$ which is infinitely generated.
Consider the subgroup $H=\langle v, w\rangle $. For any sequence
$\alpha =(\alpha _1, \ldots , \alpha _p)$, $p\in \mathbb N$, where
$\alpha _i\in \{ 0, 1\}$ for each $i=1, \ldots p$, we define an
element $t(\alpha )$ by the formula $$ t(\alpha )=w^{\alpha
_1}vw^{\alpha _2} v\ldots w^{\alpha _p}v. $$ Suppose $t(\alpha
)=t(\beta ) $ for some $\alpha =(\alpha _1, \ldots , \alpha _p)
\ne (\beta _1, \ldots , \beta _q)=\beta $. Notice that $H_{v,w}$
is normal in $H$ and $H/H_{v,w}$ is cyclic. Furthermore,
$H/H_{v,w}$ is infinite. Indeed, otherwise $H_{v,w}$ is finitely
generated). Hence $vH_{v,w}$ has infinite order when regarded as
an element of $H/H_{v,w}$. This implies $p=q$ and we have
\begin{equation}
w^{\alpha _1}vw^{\alpha _2} v\ldots w^{\alpha _p}v= w^{\beta
_1}vw^{\beta _2} v\ldots w^{\beta _p}v. \label{*}
\end{equation}
Without loss of generality, we can assume $\alpha _1\ne \beta _1$
and $\alpha _p\ne \beta _p$. Denote by $w_l$ the element
$w^{v^l}$. Then (\ref{*}) can be rewritten as $$ (w_p)^{\alpha
_1}(w_{p-1})^{\alpha _2}\ldots (w_1)^{\alpha _p}= (w_p)^{\beta
_1}(w_{p-1})^{\beta _2}\ldots (w_1)^{\beta _p}, $$ or,
equivalently, $$ (w_p)^{\alpha _1-\beta _1}=(w_{p-1})^{\beta
_2}\ldots (w_1)^{\beta _p}\left( (w_{p-1})^{\alpha _2}\ldots
(w_1)^{\alpha _p}\right) ^{-1}. $$ Note that $\alpha _1-\beta
_1=\pm 1$. Therefore,
\begin{equation}
w_p\in \langle w_1, \ldots , w_{p-1} \rangle .\label{fb}
\end{equation}
Conjugating by $v$ and using (\ref{fb}), we obtain $$
w_{p+1}=w_p^v\in \langle w_2, \ldots , w_{p} \rangle \le \langle
w_1, \ldots , w_{p-1} \rangle $$ and so on. By induction, $w_n\in
\langle w_1, \ldots , w_{p-1} \rangle $ for any $n\ge p$.
Similarly, we can show that $w_n\in \langle w_2, \ldots , w_{p}
\rangle $ for any $n\le 1$. Hence $w_n\in \langle w_1, \ldots ,
w_{p} \rangle $ for any $n\in \mathbb Z$ that contradicts to the
assumption that $H_{v,w}$ is infinitely generated.

This shows that $t(\alpha )\ne t(\beta )$ whenever $\alpha \ne
\beta $. Recall that $||v||_X\le f(s-1)$ and $||w||_X\le f(2s)$ by
Lemma 6.3. Hence we have $$ ||w^{\alpha _1}vw^{\alpha _2} v\ldots
w^{\alpha _p}v||_X\le p(||w||_X+||v||_X) \le p(f(s-1)+f(2s))\le
2pf(2s). $$ Thus, $$
\begin{array}{rl}
\gamma _H^X(n) \ge &  card\; \{ t(\alpha )\; :\; ||t(\alpha
)||_X\le n\} \\ \ge & card \left\{ (\alpha _1, \ldots , \alpha
_p)\; :\; \alpha _1, \ldots , \alpha _p \in \{ 0, 1\} , \; p\le
c(n) \right\} = 2^{c(n)},
\end{array}
$$ where $$ c(n)=\left[ \frac{n}{2f(2s)}\right].$$ Here $[x]$
means the integral part of $x$. This implies
\begin{equation}
\omega (H,X)\ge \sqrt[\beta ]{2}\label{omega}
\end{equation}
for $\beta =2f(2s)$. Note that $2f(2s)= 2(3\cdot 2^{2s-1}-2)\le
6\cdot 2^{2s-1}= 6\cdot 2^{2d+1}= 12 \cdot 4^d$. Since
(\ref{omega}) is true for an arbitrary generating set $X$, we
obtain (\ref{24}). $\Box $


\section{Applications}


In this section we consider two applications of the results
obtained above. Our goal is to prove Theorems 2.3 and 2.4.

{\bf Definition 7.1.} The {\it Cayley graph} $\Gamma = \Gamma
(G,S)$ of a group $G$ generated by a set $S$, is an oriented
labeled 1--complex with the vertex set $V(\Gamma )=G$ and the edge
set $E(\Gamma )=G\times S$. An edge $e=(g,s)\in E(\Gamma )$ goes
from the vertex $g$ to the vertex $gs$ and has the label $\phi
(e)=s$. As usual, we denote the origin and the terminus of the
edge $e$, i.e., the vertices $g$ and $gs$, by $\alpha (e)$ and
$\omega (e)$ respectively. One can endow the group $G$ (and,
therefore, the vertex set of $\Gamma $) with a {\it length
function} by assuming $\|g\|_S$, the length of an element $g\in
G$, to be equal to the length of a shortest word in the alphabet
$S\cup S^{-1}$ representing $g$.

{\bf Definition 7.2.} Two groups $G$ and $H$ with the same finite
sets of generators are called $n$--isomorphic if their Cayley
graphs restricted to balls of radius $n$ centered at the identity
are isomorphic (as labeled oriented 1--complexes). This notion
gives rise to the topology of local isomorphism (the Grigorchuk's
topology) on the set of all groups having the same collection of
generators, with subsets of $n$--isomorphic groups as the base of
neighborhoods. The local topology has interesting applications in
group theory (see \cite{Ste1}, \cite{Ste2}, \cite{Gri-84},
\cite{OsinKhyp}).

{\bf Lemma 7.1.} {\it Suppose that the group $G$ is a limit of
groups $G_i$, $i\in \mathbb N$, with respect to the Grigorchuk
topology (in this case we say that $G$ is approximated by $G_i$).
Let $X$ denote the common generating set of $G$ and $G_i$'s. Then
$$\lim\limits_{i\to\infty }\omega (G_i, X) \le \omega (G, X) .$$}

{\it Proof.} In case the group $G$ is of intermediate growth, the
proof can be found in \cite{GrHlimit}. We reproduce it here with
small changes. It is well known that $$\omega (G,X)=\inf\limits_n
\sqrt[n]{\gamma _G^X(n)}$$ (see, for example, \cite{GrH}). For any
$\varepsilon
>0$, there exists $m\in \mathbb N$ such that $$\sqrt[m]{\gamma
_G^X(m)}< \omega (G, X)+\varepsilon .$$ Further, by the definition
of the local topology, there exists $N\in \mathbb N$ such that for
every $i>N$, the balls of radius $m$ centered at the identity in
the groups $G$ and $G_i$ coincide. In particular, this implies
$$\gamma _{G_i}^X(m)=\gamma _G^X(m).$$ Finally, we have
\begin{equation}
\omega (G_i, X)=\inf\limits_n \sqrt[n]{\gamma _{G_i}^X(n)}\le
\sqrt[m]{\gamma _{G_i}^X(m)}=\sqrt[m]{\gamma _{G}^X(m)} < \omega
(G, X)+\varepsilon . \label{wGi}
\end{equation}
Since we obtain (\ref{wGi}) for an arbitrary $\varepsilon >0$, the
lemma is proved. $\Box $

{\it Proof of Theorem 2.3.} In \cite{Gri-84}, Grigorchuk constructed
a set of $3$--generated groups, $\mathcal X$, which contains a number
of groups of subexponential growth and a subset $\mathcal S$ of
elementary amenable (moreover, virtually solvable) groups of
exponential growth such that $\mathcal S$ is dense in $\mathcal X$.
Thus there exists a finitely generated group $G$ of subexponential
growth generated by a finite set $X$ and a sequence of elementary
amenable groups $G_i$ of exponential growth (generated by the same
set) that converges to $G$. By Lemma 7.1, we have $$
\lim\limits_{i\to \infty } \omega (G_i)\le \lim\limits_{i\to \infty }
\omega (G_i, X )\le \omega (G,X)=1.$$  $\Box $

{\it Proof of Theorem 2.4.} Let $G$ be a finitely generated
metabelian group. If the derived subgroup $[G,G]$ is finitely
generated, then $G$ is polycyclic. Thus we can apply Proposition
3.3 for $K=[G,G]$. In this case the group $N=G/[G, G]$ is abelian.
Substituting $d=1$ in (\ref{24}), we obtain $$ \omega (G)\ge
\sqrt[48]{2}.$$ $\Box $


\end{document}